\documentclass[a4paper,11pt]{nseJournal}
\usepackage{amsmath,amssymb,amsthm,amsfonts,stmaryrd}
\usepackage{soul}
\usepackage{array}
\usepackage{empheq}
\usepackage{enumerate}
\usepackage{enumitem}
\usepackage{graphics,graphicx}
\usepackage{longtable}
\usepackage{xcolor}
\usepackage{xfrac}
\usepackage{minibox}
\usepackage[belowskip=-8pt,font={small,it}]{caption}
\usepackage[sectionbib]{chapterbib}
\usepackage{arydshln}
\usepackage{subcaption}
\usepackage[makeroom]{cancel}

\usepackage{makecell}
\usepackage{hhline}
\usepackage{cellspace}
\usepackage{mathtools}
\DeclarePairedDelimiterX\set[1]\{\}{\nonscript\,#1\nonscript\,}
\newtheorem{remark}{Remark}
\newtheorem{algorithm}{Algorithm}

\newcommand{\tcb}{\textcolor{black}}

\begin{document}

\title{Diffusion synthetic acceleration for heterogeneous\\domains, compatible with voids}

\addAuthor{\correspondingAuthor{B. S. Southworth}}{a}
\addAuthor{Milan Holec}{b}
\addAuthor{Terry S. Haut}{b}
\correspondingEmail{ben.s.southworth@gmail.com}
\addAffiliation{a}{Department of Applied Mathematics, University of Colorado\\ Engineering Center, ECOT 225, 526 UCB, Boulder, CO 80309-0526}
\addAffiliation{b}{Lawrence Livermore National Laboratory,\\7000 East Ave, Livermore, CA 94550.}

\addKeyword{transport}
\addKeyword{discrete ordinates}
\addKeyword{high-order}
\addKeyword{sweep}
\addKeyword{unstructured}

\author{
 \vspace{10mm} \\
   B. S. Southworth$^a$,
   Milan Holec$^b$,
   T. S. Haut$^b$
   \\
   \\
   \normalsize{
    \begin{itshape}
    \begin{tabular}[t]{c}
           $^a$ Department of Applied Mathematics, University of Colorado\\ Engineering Center, ECOT 225, 526 UCB, Boulder, CO 80309-0526\\
           $^b$ Lawrence Livermore National Laboratory,\\7000 East Ave, Livermore, CA 94550.
    \end{tabular}
    \end{itshape}
    }
     \vspace{10mm}
    \\
    $^*$Email: \url{mailto:ben.s.southworth@gmail.com}
    }

 \date{
        \vspace{15mm}
        \begin{tabular}[t]{lc}
        Number of pages:& \pageref*{LastPage} \\
        Number of tables:& \totaltables \\
        Number of figures:& \totalfigures
        \end{tabular}
    }
 \clearpage\maketitle
    \thispagestyle{empty}

\begin{abstract}
A standard approach to solving the S$_N$ transport equations is to use source iteration with diffusion synthetic acceleration
(DSA). Although this approach is widely used and effective on many problems, there remain some practical issues
with DSA preconditioning, particularly on highly heterogeneous domains. For large-scale parallel simulation, it is
critical that both (i) preconditioned source iteration converges rapidly, and (ii) the action of the DSA preconditioner
can be applied using fast, scalable solvers, such as algebraic multigrid (AMG). For heterogeneous domains, these
two interests can be at odds. In particular, there exist DSA diffusion discretizations that can be solved rapidly using
AMG, but they do not always yield robust/fast convergence of the larger source iteration. Conversely, there exist
robust DSA discretizations where source iteration converges rapidly on difficult heterogeneous problems, but fast
parallel solvers like AMG tend to struggle applying the action of such operators. Moreover, very few current methods
for the solution of deterministic transport are compatible with voids. This paper develops a new heterogeneous DSA
preconditioner based on only preconditioning the optically thick subdomains. The resulting method proves robust
on a variety of heterogeneous transport problems, including a linearized hohlraum mesh related to inertial
confinement fusion. Moreover, the action of the preconditioner is easily computed using $\mathcal{O}(1)$ AMG
iterations, {convergence of the transport iteration typically requires $2-5\times$ less iterations than
current state-of-the-art ``full DSA,'' and the proposed method is} trivially compatible
with voids. On the hohlraum problem, rapid convergence is obtained by preconditioning less than 3\% of
the mesh elements with $5-10$ AMG iterations. 
\end{abstract}

\section{Introduction} \label{sec:intro}

Solving the particle transport equations arises in numerous fields of research, such as nuclear reactor design,
inertial confinement fusion (ICF), and medical imaging. Despite decades of research, however, it remains
a computationally expensive and challenging problem. The simplified steady state, mono-energetic neutral-particle transport
equation for the spatially and angularly dependent angular flux, $\psi$, is given by
\begin{align}\label{eq:eq}
\boldsymbol{\Omega} \cdot \nabla{\psi}(\mathbf{x},\boldsymbol{\Omega}) + \sigma_t(\mathbf{x}) {\psi}(\mathbf{x},\boldsymbol{\Omega})
	=  \frac{\sigma_s(\mathbf{x})}{4\pi} \int_{\mathbb{S}^2} {\psi}(\mathbf{x},\boldsymbol{\Omega'})
	 d\boldsymbol{\Omega'} + q(\mathbf{x},\boldsymbol{\Omega}).
\end{align}
Here, $0\leq \sigma_s\leq \sigma_t$ are the spatially dependent scattering and total cross sections, respectively, and
we assume isotropic scattering for simplicity.
This equation is fundamental to models of neutral particle transport, such as neutron or photon transport, and
also a prerequisite for solving more complicated physical models such as thermal radiative transfer. 
The left-hand side of \eqref{eq:eq} is a 3D advection-reaction equation in direction $\boldsymbol{\Omega}$
and the right-hand side an integral-operator coupling over direction. {More complex physical models
introduce a nonlinearity in temperature on the right-hand side, while including time and energy leads to seven
dimensions, overall requiring massively parallel simulations.}

In this paper, we consider an S$_N$ angular discretization of the integral operator in \eqref{eq:eq}.
{Using $n$ directions in the angular quadrature set} yields the following semi-discrete form of \eqref{eq:eq},
\begin{align}
\begin{split}
\boldsymbol{\Omega}_{d} \cdot \nabla_{\mathbf{x}}\psi_{d}(\mathbf{x}) + \sigma_{t}\left(\mathbf{x}\right)\psi_{d}(\mathbf{x}) 
	& = q_{d}(\mathbf{x}) + \frac{\sigma_s(\mathbf{x})}{4\pi} \sum_{d'=1}^{n}\omega_{d'}\psi_{d'}\left(\mathbf{x}\right) , 
	\,\,\,\, \mathbf{x} \in \mathcal{D}, \\
\psi_{d}(\mathbf{x}) & =\psi_{d,\text{inc}}(\mathbf{x}),\,\,\,\,\mathbf{x}\in \partial{\mathcal{D}} \,\,\,\,\text{and}\,\,\,\,
	\mathbf{n}(\mathbf{x})\cdot\boldsymbol{\Omega}_{d}<0,
\label{eq:sntransport}
\end{split}
\end{align}
where subscript $d$ denotes a fixed angle in the $S_N$ angular discretization, with weight $\omega_d$ and
direction $\boldsymbol{\Omega}_{d}$. Here $\psi_{d}(\mathbf{x})$ denotes the angular flux associated with
direction $d$, and the scalar flux is then defined as 
\begin{align*}
\varphi\left(\mathbf{x}\right) := \sum_{d'=1}^{n}w_{d'}\psi_{d'}\left(\mathbf{x}\right).
\end{align*}
Due to the high-dimensionality of the transport equations, memory is one of the fundamental constraints in
solving them numerically. Fortunately, it is generally the case that $\{\psi_d\}$ can be eliminated from the 
problem, and the scalar flux $\varphi$ iterated to convergence, thus only requiring storage of one vector
on the spatial domain rather than one for each angle. This is discussed formally in a linear algebraic 
setting in Section \ref{sec:alg:si}.

A common approach to solve \eqref{eq:eq} is based on a fixed-point iteration, where each iteration 
inverts the left-hand side for each direction in the S$_N$ angular discretization. 
{We will refer to this set of inversions as a transport update, which is the dominant computational
cost of solving the transport equations. Upwinded discontinuous discretizations are particularly
amenable to an efficient update (and commonly used for this reason), because the linear system
for each direction can be solved directly using a forward solve. A transport \textit{sweep}
then consists of a parallel implementation of simultaneous forward solves for all directions
in the quadrature set.} For $\sigma_t \ll 1$ or $\sigma_s\ll\sigma_t$, {fixed-point iteration based on
a transport update converges rapidly.} For
$\sigma_t \approx\sigma_s\gg 1$, the fixed-point operator is very slow to converge, but can be effectively
preconditioned with so-called diffusion synthetic acceleration (DSA). DSA corresponds to using an appropriate discrete diffusion operator
as a preconditioner for the fixed-point iteration on the scalar flux, which is the focus of this paper. 

DSA was developed in separate works in the 1960s as an acceleration technique for solving the
transport equations by approximating the behavior of the scalar flux in optically thick materials with
a diffusion equation \cite{kopp1963synthetic,lebedev1964iterative}. DSA was analyzed and improved for specific
discretizations in the 1970s (for example, \cite{Alcouffe:1977wv}) and a new scaling of the transport
equations was introduced in \cite{Larsen:1982hh} to analyze DSA preconditioning. Since, DSA
has seen significant research over the years. Some of the most important, general observations
include (i) the need for a so-called consistency between the discretization of the linear transport equation
and the diffusion preconditioning, and (ii) the degraded effectiveness of DSA in heterogeneous
media \cite{Warsa:2004gt}. In particular, DSA is known to be effective and easy to apply on
homogeneous and optically thick domains, $\sigma_t \gg 1$. However, it is not uncommon for
$\sigma_s$ and $\sigma_t$ to vary by many orders of magnitude over a spatial domain, and
such heterogeneous domains remain difficult for existing methods.

The practical difficulties of DSA on heterogeneous domains come from trying to satisfy both of
the requirements for a fast, parallel transport simulation:
\begin{enumerate}
\itemsep-0em
\item The DSA discretization must be an effective preconditioner of the S$_N$ transport 
equations, resulting in rapid convergence.
\item The DSA preconditioning must be able to be applied in a fast, parallel manner.
\end{enumerate}
{In most cases these two goals can be satisfied individually, and for some problems
they can both be
satisfied.} However, it is often difficult or impossible to satisfy both using existing
techniques on highly heterogeneous domains. Part of this problem comes from the
requirement of consistency between the transport and diffusion discretizations. In particular,
the transport equations are often discretized using some form of upwind
discontinuous discretization. It turns out the compatible DSA preconditioner involves inverting a symmetric interior
penalty type diffusion discretization.
It is well known though that fast solvers such as algebraic multigrid (AMG) tend to struggle
on discontinuous discretizations, even of elliptic problems. This is compounded by the
fact that these compatible DSA diffusion discretizations can have very poor
conditioning \cite{18dsa} {(worse than a standard diffusion conditioning of
$\sim\mathcal{O}(1/h^2)$).}

Thus the dichotomy is as follows: rapid convergence on highly heterogeneous domains
requires DSA based on specific discontinuous diffusion discretizations. However, such
discretizations are between difficult and just not amenable to existing fast, parallel linear
solvers, such as AMG, particularly if using existing parallel linear solver libraries such as
\textit{hypre} \cite{Falgout:2002vu}. A number of works have considered how to construct
more effective algebraic solvers for such diffusion discretizations, including in the 
transport community \cite{Warsa:2003ug,OMalley:2017ck} as well as the linear
solver community \cite{Olson:2011ju,Bastian:2012hk,Antonietti:2016fa}, but 
to our knowledge (and experience testing various methods) heterogeneous problems remain
difficult for existing methods, particularly methods that are publicly available and implemented
in parallel. Furthermore, it is generally the case that DSA bilinear forms have a 
diffusion coefficient proportional to $1/\sigma_t(\mathbf{x})$, which is not well-defined if there
are any regions of vacuum, $\sigma_t = 0$. Some works have developed 
appropriate linear and nonlinear DSA modifications when using the self-adjoint angular
flux formulation of transport \cite{Wang:2017jx,Hammer:2019fj}, but to our knowledge,
DSA applied to the S$_N$ transport equations in domains with voids remains an open
question. 

This paper derives a new DSA-like preconditioner (not discretization) for heterogeneous
domains, which is (i) easier to apply in a fast and scalable manner with existing solvers
than standard DSA, (ii) a better preconditioner than even the most robust standard DSA
preconditioners we have tested, and (iii) amenable to vacuum. Conceptually, the new
algorithm is quite simple, and more-or-less corresponds to only applying DSA preconditioning
on ``thick'' regions in the domain.\footnote{Recent work in \cite{19char} has seen success applying
similar ideas to precondition the radiation diffusion equations by only solving the diffusion
discretization on a physical subdomain.} Using variations in a crooked pipe benchmark problem,
the new heterogeneous DSA method reduces the total number of iterations to convergence
by $5-6\times$ for some problems. Moreover, the heterogeneous DSA matrices are more tractable
to solve using AMG, in a number of tests converging in $\mathcal{O}(1)$ iterations when AMG
was unable to solve the full DSA matrix. The new method is relatively non-intrusive, in the sense
that it can be added to existing libraries that support DSA with minimal work. Conceptually, the
preconditioner is independent of the discretization. However, a discussion on implementation and
numerical results are presented based on a discontinuous Galerkin (DG) discretization of the
S$_N$ transport equations and corresponding DG DSA discretization \cite{18dsa,Wang:2010dva}. 

Section \ref{sec:alg} presents standard source iteration and DSA in the context of block preconditioners. This motivates
a similar analysis applied specifically to DSA preconditioning, and the development of the new heterogeneous DSA
preconditioning in Section \ref{sec:alg:het}. An implementation involving the \textit{hypre} library and DG spatial
discretizations is discussed in Section \ref{sec:imp}, including a generalization of the modified penalty coefficient 
introduced in \cite{Wang:2010dva} for high-order curvilinear meshes. 
Numerical results then demonstrate the new method on heterogeneous domains in Section \ref{sec:results},
including on variations of a crooked pipe problem \cite{gentile2001implicit,SmedleyStevenson:2015wa}, and 
a linearized hohlraum capsule. For the hohlraum, in particular, rapid convergence is obtained when applying
(heterogeneous) DSA preconditioning to less than $ 3\%$ of mesh elements! Some conclusions are given in
Section \ref{sec:conc}.

\section{Preconditioning linear transport} \label{sec:alg}

\subsection{Review of $2\times 2$ block preconditioners}\label{sec:alg:2x2}

Source iteration and diffusion synthetic acceleration (DSA) applied to linear $S_N$ transport can be
seen in a linear-algebraic framework as a block preconditioning of a $2\times 2$ block operator. A
linear-algebraic perspective on transport iterations was introduced as early as 1989 in \cite{Faber:1989wo}
\footnote{Unfortunately, the only available copy of \cite{Faber:1989wo} appears to be a low-resolution,
scanned in copy, which is difficult to read.}, and transport iterations are often expressed in algebraic 
operator form (for example, see \cite{Wang:2010dva}). The specific context of $2\times 2$ block
preconditioners is not something seen often in the literature, but provides valuable insight on the
preconditioning of heterogeneous domains.

As a background, consider a $2\times 2$ block matrix $A$ and lower-triangular preconditioner $L$,
\begin{align}\label{eq:block}
A & = \begin{bmatrix} A_{11} &A_{12} \\ A_{21} & A_{22}\end{bmatrix}, \hspace{8ex}
	L = \begin{bmatrix} A_{11} & \mathbf{0} \\ A_{21} & \widehat{\mathcal{S}} \end{bmatrix},
\end{align}
where $\mathcal{S} := A_{22} - A_{21}A_{11}^{-1}A_{12}$ is the Schur complement of $A$ in the
$(2,2)$-block, and $\widehat{\mathcal{S}}$ some approximation to $\mathcal{S}$. In the
simplest case of a block Gauss-Seidel like iteration, we have $\widehat{\mathcal{S}} := A_{22}$.
It can be shown that fixed-point or Krylov iterations applied to $A\mathbf{x} = \mathbf{b}$ with preconditioner
$L^{-1}$ converge to some tolerance $< C\epsilon$ after $k$ iterations {if and only if} equivalent iterations
applied to a Schur complement problem $\mathcal{S}\mathbf{x}_c = \mathbf{b}_c$, with preconditioner
$\widehat{\mathcal{S}}^{-1}$, converge to tolerance $\epsilon$ after $k$ iterations, where 
$C \sim\mathcal{O}(\|A_{11}^{-1}A_{12}\|)$ \cite{2x2}.  As it turns out, source iteration and DSA
preconditioning correspond exactly to a block lower-triangular preconditioning, as in \eqref{eq:block},
and considering transport iterations in the context of block preconditioners allows us to
derive a natural approach to DSA preconditioning in heterogeneous media.

\subsection{Source iteration as a block preconditioner}\label{sec:alg:si}

Now suppose the left-hand side of \eqref{eq:sntransport} is discretized in space for each direction $d$, with
corresponding discrete spatial operator $\mathcal{L}_d \sim \boldsymbol{\Omega_d}\cdot\nabla_{\mathbf{x}} + \sigma_t$,
and let $\Sigma_s$ and $\Sigma_t$ denote mass matrices with respect to coefficients $\sigma_s(\mathbf{x})/4\pi$
and $\sigma_t(\mathbf{x})$, respectively. Further, let bold $\boldsymbol{\psi}_d$ and $\boldsymbol{q}_d$
denote the discrete vector representations of $\psi_{d}(\mathbf{x})$ and $q_d(\mathbf{x})$. Then, the full
discretized set of equations can be written as a block linear system,
\begin{align} \label{eq:system}
 \left[\begin{array}{@{}ccc : c@{}}
	\mathcal{L}_1 &&& - \Sigma_s \\
	&\ddots&& \vdots \\
	&&\mathcal{L}_n& -\Sigma_s \\\hdashline
	-\omega_1 I &  \hdots &  -\omega_n I &  I \end{array}\right]
\begin{bmatrix} \boldsymbol{\psi}_1 \\ \vdots \\ \boldsymbol{\psi}_n \\\hdashline \boldsymbol{\varphi} \end{bmatrix} =
	\begin{bmatrix}\mathbf{q}_1 \\ \vdots \\ \mathbf{q}_n \\ \hdashline\mathbf{0} \end{bmatrix}.
\end{align}
The dotted lines indicate how the $S_N$ transport equations can be expressed as a $2\times 2$ block array
\eqref{eq:block}: $A_{11}$ is the block-diagonal set of spatial operators, $\{\mathcal{L}_d\}$, $A_{22} = I$,
and the off-diagonal operators $A_{12}$ and $A_{21}$ correspond to scattering and quadrature weights,
respectively. 

One of the standard approaches to solve \eqref{eq:system} in transport simulations is to update the scalar flux
based on the current angular flux,
\begin{align}\label{eq:psi}
\boldsymbol{\psi}_d^{(i+1)} & = \mathcal{L}_d^{-1}(\mathbf{q}_d + \Sigma_s\boldsymbol{\varphi}^{(i)}),
\end{align}
for all directions $d=1,...,n$. With updated scalar flux for each direction, the angular flux is then updated by
summing the scalar flux over quadrature weights $\{\omega_d\}$,
\begin{align}\label{eq:phi}
\boldsymbol{\varphi}^{(i+1)} & = \sum_{d=1}^n \omega_d \boldsymbol{\psi}_d^{(i+1)} =
	\sum_{d=1}^n \omega_d\mathcal{L}_d^{-1}\left(\mathbf{q}_d + \Sigma_s\boldsymbol{\varphi}^{(i)}\right).
\end{align}
This process is repeated and is the classical ``source iteration.'' Algebraically, source iteration is exactly
a fixed-point iteration with block lower-triangular preconditioning:
{\small
\begin{align*}
\begin{bmatrix} \boldsymbol{\psi}_1^{(i+1)} \\ \vdots \\ \boldsymbol{\psi}_n^{(i+1)} \\ \boldsymbol{\varphi}^{(i+1)} \end{bmatrix}
& = \begin{bmatrix} \boldsymbol{\psi}_1^{(i)} \\ \vdots \\ \boldsymbol{\psi}_n^{(i)} \\ \boldsymbol{\varphi}^{(i)} \end{bmatrix} + 
	\begin{bmatrix} \mathcal{L}_1 &&&  \\ &\ddots&&  \\ &&\mathcal{L}_n& \\
	-\omega_1 I &  \hdots &  -\omega_nI &  I \end{bmatrix}^{-1} \left( 
	\begin{bmatrix} \mathbf{q}_1 \\ \vdots \\ \mathbf{q}_n \\ \mathbf{0}\end{bmatrix} -
	\begin{bmatrix} \mathcal{L}_1 &&& - \Sigma_s \\ &\ddots&& \vdots \\ &&\mathcal{L}_n& -\Sigma_s \\
		-\omega_1 I &  \hdots &  -\omega_n I &  I \end{bmatrix}
	\begin{bmatrix} \boldsymbol{\psi}_1^{(i)} \\ \vdots \\ \boldsymbol{\psi}_n^{(i)} \\ \boldsymbol{\varphi}^{(i)} \end{bmatrix} \right).
\end{align*}
}
Expanding, one arrives at exactly the two stage iteration for $\{\boldsymbol{\psi}_d\}$ and $\boldsymbol{\varphi}$
introduced in \eqref{eq:psi} and \eqref{eq:phi}. In this form, it is also clear how Krylov methods can be applied to
accelerate convergence of source iteration \cite{warsa2003improving,Warsa:2004gt}. In particular, Krylov methods
require computing the action of the operator, $A$, and preconditioner, $M^{-1}$, on vectors, where
\begin{align*}
A = \begin{bmatrix} \mathcal{L}_1 &&& - \Sigma_s \\ &\ddots&& \vdots \\ &&\mathcal{L}_n& -\Sigma_s \\
		-\omega_1 I &  \hdots &  -\omega_n I &  I \end{bmatrix}, 
\hspace{5ex}
M^{-1} = \begin{bmatrix} \mathcal{L}_1^{-1} &&&  \\ &\ddots&&  \\ &&\mathcal{L}_n^{-1}& \\
	\omega_1 \mathcal{L}_1^{-1} &  \hdots &  \omega_n\mathcal{L}_n^{-1} & I \end{bmatrix} .
\end{align*}

Note that in source iteration, the transport update inverts $\mathcal{L}_d^{-1}$ for all $d=1,...,n$,
meaning that in the $2\times 2$ block sense, $A_{11}$ is inverted exactly. Recall from Section
\ref{sec:alg:2x2} that convergence of fixed-point or Krylov iterations with a $2\times 2$ block lower-triangular
preconditioner are then defined by equivalent iterations on the preconditioned Schur complement problem.
In this case, in the notation of \eqref{eq:block}, we have $\widehat{\mathcal{S}} = A_{22} = I$, and
the preconditioned Schur complement takes the form $\mathcal{S}\boldsymbol{\varphi} = \mathbf{b}$, where
\begin{align}
\mathcal{S} & := I - \begin{bmatrix} \omega_1 I &  \hdots &  \omega_n I \end{bmatrix}
	\begin{bmatrix}  \mathcal{L}_1^{-1} & \\ & \ddots &\\ &&\mathcal{L}_n^{-1} \end{bmatrix} \begin{bmatrix}
	\Sigma_s \\  \vdots \\  \Sigma_s\end{bmatrix}
	= I - \sum_{d=1}^n \omega_d\mathcal{L}_d^{-1}\Sigma_s,\label{eq:schur}
\end{align}
and $\mathbf{b} := \sum_{d=1}^n \omega_d\mathcal{L}_d^{-1}\mathbf{q}_d.$
If $\mathcal{S}$ is well conditioned, such as in the optically thin case when $\Sigma_s \ll 1$, we can solve 
$\mathcal{S}\boldsymbol{\varphi} = \mathbf{b}$ via a simple Richardson iteration,
\begin{align}
\boldsymbol{\varphi}^{(i+1)} & = \boldsymbol{\varphi}^{(i)} + \mathbf{b} - \mathcal{S}\boldsymbol{\varphi}^{(i)} 
	 =    \sum_{d=1}^n \omega_d\mathcal{L}_d^{-1}\left(\mathbf{q}_d + \Sigma_s\boldsymbol{\varphi}^{(i)}\right). \label{eq:source0}
\end{align}

If $\mathcal{S}$ is ill conditioned, such as in optically thick or heterogeneous media, some form of
preconditioner for $\mathcal{S}$ is necessary for fast convergence. In the larger system, this simply
replaces the $I$ in the lower right $(2,2)$-block of the preconditioner, $M^{-1}$, with some approximation
to $\mathcal{S}$. Denoting this approximation $\mathcal{D}$,
the preconditioner takes the form
\begin{align}
M^{-1} & = \begin{bmatrix} \mathcal{L}_1 &&&  \\
&\ddots&&  \\
&&\mathcal{L}_n& \\
-\omega_1 I &  \hdots &  -\omega_n I &  \mathcal{D} \end{bmatrix}^{-1}
= 
\begin{bmatrix} I &&&  \\
&\ddots&&  \\
&&I& \\
&  & &  \mathcal{D}^{-1} \end{bmatrix}
\begin{bmatrix} \mathcal{L}_1^{-1} &&&  \\
&\ddots&&  \\
&&\mathcal{L}_n^{-1}& \\
\omega_1 \mathcal{L}_1^{-1} &  \hdots &
	\omega_n \mathcal{L}_n^{-1} &  I \end{bmatrix}.
\label{eq:dsa_prec}
\end{align}
Now, convergence of fixed-point or Krylov iterations applied to the larger system \eqref{eq:system} is
defined by convergence of equivalent iterations applied to $\mathcal{D}^{-1}\mathcal{S}$. In assuming that 
$\mathcal{L}_d$ is inverted for all $d$ (rather than, say, an on-processor solve), $\{\boldsymbol{\psi}_d\}$
can be eliminated from the system, yielding a preconditioned fixed-point iteration on $\boldsymbol{\varphi}$,
\begin{align}
\boldsymbol{\varphi}^{(i+1)} 
& = \boldsymbol{\varphi}^{(i)} - \mathcal{D}^{-1}\left(\boldsymbol{\varphi}^{(i)} - \sum_{d=1}^n
	\omega_d\mathcal{L}_d^{-1}\left( \Sigma_s\boldsymbol{\varphi}^{(i)} + \mathbf{q}_d\right)\right). \label{eq:DSA_fp}
\end{align}
At any point, an approximation to the angular flux can be computed via equation (\ref{eq:psi})
(and is implicitly computed within every iteration). Similar to the case of source iteration, it is straightforward
to analyze and implement Krylov methods here as well, where $\mathcal{S}$ is the operator and $\mathcal{D}^{-1}$
the preconditioner. It is worth pointing out that for linear S$_N$ transport as considered here, if angular flux
$\{\boldsymbol{\psi}_d^{(i+1)}\}$ is eliminated from the system, it is clear that Krylov convergence is defined by
convergence on the Schur complement (scalar flux). {For nonlinear or time-dependent transport problems,
where the angular flux cannot be fully eliminated}, this relation is less obvious without considering
iterations in the context of $2\times 2$ block-preconditioning and appealing to the theory on block preconditioners
and Krylov methods \cite{2x2}.

{Note that in practice, $\mathcal{D}$ is not directly a diffusion operator, rather {$\mathcal{D}^{-1}
= I + D^{-1}\Sigma$ is an additive preconditioner, with diffusion operator $D$ and mass matrix
$\Sigma\sim \sigma_s(\mathbf{x})$ \tcb{or $\Sigma\sim \sigma_t(\mathbf{x})$}.\footnote{\tcb{Using
$\Sigma\sim \sigma_s(\mathbf{x})$ is more common in literature; here we use $\Sigma\sim
\sigma_t(\mathbf{x})$ as derived in \cite{18dsa}; in any case, they can be shown to have the
same asymptotic efficiency in the ``thick'' limit.}} There are various ways to 
motivate the additive preconditioner. Conceptually, a diffusion operator can be derived in which
$D^{-1}\Sigma\mathcal{S}$ is well-conditioned in the optically thick limit of mean free path $\varepsilon \ll 1$
(for DG discretizations, see \cite{Wang:2010dva,18dsa}). However, in optically thin material, 
$\Sigma_s \ll 1$, $D^{-1}\Sigma\mathcal{S}$ can be ill conditioned because $\mathcal{S}$ is already well
conditioned in thin material, but the spectrum of $D^{-1}$ goes to zero for high-frequency modes.
Adding the identity shifts such eigenmodes so that $(I+D^{-1}\Sigma)\mathcal{S}$ is well conditioned in
thin regimes, while for thick, note that if $D^{-1}\Sigma\mathcal{S}$ is positive and well-conditioned,
then $(I+D^{-1}\Sigma)\mathcal{S}$ is also positive and well-conditioned.}

\begin{remark}
If $\mathcal{L}_d$ is not inverted exactly due to, for example, {on-processor block Jacobi iterations
instead of a global inverse,} or cycle-breaking on high-order curvilinear meshes \cite{19sweep}, it is
straightforward to work out a reduced
iteration similar to above, which stores $\boldsymbol{\varphi}$ and the degrees-of-freedom (DOFs) of
$\boldsymbol{\psi}_d$ necessary to update each direction $d$. Extending this to the Krylov setting takes
a little more care, but follows in an analogous manner \cite{wang2009adaptive}.
\end{remark}

\subsection{Heterogeneous DSA with two regions}\label{sec:alg:het}

Now consider the case of heterogeneous domains, and suppose we partition the domain into ``thin'' and ``thick''
regions. In particular, let DOFs/elements in which $\sigma_s(\mathbf{x}) < \eta$ be considered thin,
denoted by subscript $s$, and DOFs/elements such that $\sigma_s(\mathbf{x}) \geq \eta$ be considered
thick, denoted with subscript $f$. Then, order all matrices $\{\mathcal{L}_d\}$ in a block ordering
\begin{align*}
\mathcal{L}_d & = \begin{bmatrix} {L}_{ss} & {L}_{sf} \\ {L}_{fs} & {L}_{ff} \end{bmatrix}.
\end{align*}
Note that these region do not have to be contiguous. Moving forward, subscripts $d$ denoting angle are
dropped on submatrices $L_{ss}, L_{sf}, L_{fs}$, and $L_{ff}$, for ease of notation. Note, however, that
summations over direction $d$ have an implied subscript on submatrices of $\mathcal{L}_d$. 

Now, recall the two Schur complements of $\mathcal{L}_d$,
\begin{align}\label{eq:schurinv}
\mathcal{S}_{ss} & = {L}_{ss} - {L}_{sf}{L}_{ff}^{-1}{L}_{fs}, \hspace{5ex}
\mathcal{S}_{ff} = {L}_{ff} - {L}_{fs}{L}_{ss}^{-1}{L}_{sf},
\end{align}
and the closed form inverse of a $2\times 2$ block matrix {(see, for example, \cite[3.2.11]{golub})},
\begin{align}
\mathcal{L}_d^{-1} & = \begin{bmatrix}  \mathcal{S}_{ss}^{-1} & -L_{ss}^{-1}L_{sf}\mathcal{S}_{ff}^{-1}  \\
	-\mathcal{S}_{ff}^{-1}L_{fs}L_{ss}^{-1} & \mathcal{S}_{ff}^{-1} \end{bmatrix}, \label{eq:binv}
\end{align}
where
\begin{align*}
\mathcal{S}_{ss}^{-1} & = L_{ss}^{-1} + L_{ss}^{-1}L_{sf}\mathcal{S}_{ff}^{-1}L_{fs}L_{ss}^{-1}, \hspace{5ex}
\mathcal{S}_{ff}^{-1}  = L_{ff}^{-1} + L_{ff}^{-1}L_{fs}\mathcal{S}_{ss}^{-1}L_{sf}L_{ff}^{-1}
\end{align*}
{are interdependent expressions for the inverse of a Schur complement.} Using \eqref{eq:binv}, we can expand the Schur complement \eqref{eq:schur} in block form as
\begin{align}
\mathcal{S} & = I - \sum_{d=1}^n \omega_d\mathcal{L}_d^{-1}\Sigma_s \nonumber\\
& = I - \sum_{d=1}^n \omega_d
	 \begin{bmatrix}  \mathcal{S}_{ss}^{-1} & -L_{ss}^{-1}L_{sf}\mathcal{S}_{ff}^{-1}  \\
	-\mathcal{S}_{ff}^{-1}L_{fs}L_{ss}^{-1} & \mathcal{S}_{ff}^{-1} \end{bmatrix}
	\begin{bmatrix} \Sigma_{s,s} \\ & \Sigma_{s,f}\end{bmatrix}.\label{eq:s_block}
\end{align}
Observe that $\Sigma_{s}$ is a column scaling. Suppose the thin region is actually vacuum, in which
case $\Sigma_{s,s} = \mathbf{0}$. Then,
\begin{align}
\mathcal{S} & = I - \sum_{d=1}^n \omega_d
	 \begin{bmatrix}  \mathbf{0} & -L_{ss}^{-1}L_{sf}\mathcal{S}_{ff}^{-1}\Sigma_{s,f}  \\
	\mathbf{0} & \mathcal{S}_{ff}^{-1}\Sigma_{s,f} \end{bmatrix} \nonumber\\
& =  \begin{bmatrix} I & \sum_{d=1}^n \omega_dL_{ss}^{-1}L_{sf}\mathcal{S}_{ff}^{-1}\Sigma_{s,f}  \\
	\mathbf{0} & I - \sum_{d=1}^n \omega_d\mathcal{S}_{ff}^{-1}\Sigma_{s,f} \end{bmatrix}. \label{eq:Stri}
\end{align}

Looking at \eqref{eq:Stri}, a natural choice of preconditioner for $\mathcal{S}$ is a block upper-triangular
preconditioner, with ``blocks'' indicating thin and thick regions of the domain. Now we only need to
precondition the lower right (thick) block of \eqref{eq:Stri}, $ I - \sum_{d=1}^n \omega_d\mathcal{S}_{ff}^{-1}\Sigma_{s,f}$.
{A natural choice here is to approximate $\mathcal{S}_{ff}^{-1}$ in the lower diagonal block
(but not in the off-diagonal block) by $L_{ff}^{-1}$, in some sense a first-order approximation
from \eqref{eq:schurinv}.} This yields,
\begin{align}\label{eq:approxS0}
\mathcal{S} & \approx  \begin{bmatrix} I & \sum_{d=1}^n \omega_dL_{ss}^{-1}L_{sf}\mathcal{S}_{ff}^{-1}\Sigma_{s,f}  \\
	\mathbf{0} & I - \sum_{d=1}^n \omega_d{L}_{ff}^{-1}\Sigma_{s,f} \end{bmatrix}.
\end{align}
Now, disregarding boundary conditions, $I - \sum_{d=1}^n \omega_d{L}_{ff}^{-1}\Sigma_{s,f}$ corresponds
to $S_N$ transport only on the thick region, for which we know that diffusion is an effective preconditioning.
Letting $D_{ff}^{-1}$ denote an appropriate diffusion discretization over the thick region,
we use the preconditioner $I + D_{ff}^{-1}\Sigma_{t,f} \approx (I - \sum_{d=1}^n \omega_d{L}_{ff}^{-1}\Sigma_{s,f})^{-1}
\approx (I - \sum_{d=1}^n \omega_d\mathcal{S}_{ff}^{-1}\Sigma_{s,f})^{-1}$,
which can be applied as the first step in a block upper triangular preconditioning,
\begin{align}
\mathcal{D}_T^{-1} & = \begin{bmatrix} I & \sum_{d=1}^n \omega_dL_{ss}^{-1}L_{sf}\mathcal{S}_{ff}^{-1}\Sigma_{s,f}  \\
	\mathbf{0} & (I + D_{ff}^{-1}\Sigma_{t,f})^{-1} \end{bmatrix}^{-1} \nonumber\\
& = \begin{bmatrix} I & -\sum_{d=1}^n \omega_dL_{ss}^{-1}L_{sf}\mathcal{S}_{ff}^{-1}\Sigma_{s,f} \\
	\mathbf{0} & I \end{bmatrix} \begin{bmatrix} I & \mathbf{0} \\ \mathbf{0} & I + D_{ff}^{-1}\Sigma_{t,f}\end{bmatrix} \label{eq:approxS}.
\end{align}
The resulting preconditioning corresponds to a DSA preconditioning only on the thick region (including interface),
followed by using the updated thick solution as an additive correction to the thin region. Note that applying DSA
only on the thick region is \textit{not} the same as inverting a global diffusion operator and only updating the thick DOFs. 

To apply the additive correction in the left block of \eqref{eq:approxS}, following DSA on the thick region,
note from \eqref{eq:binv} that 
\begin{align*}
\begin{bmatrix} I & -\sum_{d=1}^n \omega_dL_{ss}^{-1}L_{sf}\mathcal{S}_{ff}^{-1}\Sigma_{s,f} \\
	\mathbf{0} & I \end{bmatrix}
& = I + \begin{bmatrix} I & \mathbf{0} \\ \mathbf{0} & \mathbf{0} \end{bmatrix}
	\sum_{d=1}^n \omega_d\mathcal{L}_d^{-1}\Sigma_s
	\begin{bmatrix} \mathbf{0} & \mathbf{0} \\ \mathbf{0} & I \end{bmatrix}.
\end{align*}
This is nothing more than a modified transport update with zero right-hand side. For all angles, the
angular flux is first set to zero on the thin region, then a transport update is applied, inverting $\mathcal{L}_d$
for each direction $d$, but only accumulating the solution corrections in thin DOFs. 
The action of the preconditioner \eqref{eq:approxS} can be expressed as two steps:\\
\begin{algorithm}[Triangular heterogeneous DSA]\label{alg:tri}
\text{ }
\begin{enumerate}
\item Let $\boldsymbol{\varphi}_f$ denote $\boldsymbol{\varphi}$ restricted to the thick region, and
update $\boldsymbol{\varphi}_f \mapsfrom \boldsymbol{\varphi}_f + D_{ff}^{-1}\Sigma_{t,f}\boldsymbol{\varphi}_f$,
where $D_{ff}$ is a DSA diffusion preconditioner over the thick region.

\item Define $\overline{\boldsymbol{\varphi}} := \begin{bmatrix}\mathbf{0} \\ {\boldsymbol{\varphi}}_{f}\end{bmatrix}$,
and let $\boldsymbol{\varphi}_s$ denote $\boldsymbol{\varphi}$ restricted to the thin region. For each direction $d$,
update $ \boldsymbol{\varphi}_s$ via
\begin{align*}
\boldsymbol{\varphi}_s \mapsfrom \boldsymbol{\varphi}_s +
	\left[\omega_d\mathcal{L}_d^{-1}\Sigma_s\overline{\boldsymbol{\varphi}} \right]_s.
\end{align*}
\end{enumerate}
\end{algorithm}

In the case of $\Sigma_{s,s} = \mathbf{0}$, $\mathcal{S}$ is indeed block upper triangular, and the
convergence of block-triangular preconditioned minimal-residual methods is defined by the approximation
of the $(2,2)$-block in $\mathcal{S}$. For $\Sigma_{s,s} > 0$, $\mathcal{S}$ is no longer triangular, and
convergence is then defined by preconditioning of the (2,2)-Schur complement
of $\mathcal{S}$ \cite{2x2}. However, when $\Sigma_{s,s}$ is small, this preconditioning will likely be
comparable to the preconditioning of the (2,2)-block in $\mathcal{S}$, analogous to when $\Sigma_{s,s}
= \mathbf{0}$. That is, convergence of heterogeneous-DSA-preconditioned source iteration is
expected to be comparable for general $\Sigma_{s,s} \ll1$, which is consistent with numerical results
in Section \ref{sec:results}. 

In many cases block-triangular preconditioners offer faster convergence than block-diagonal, with marginal
additional cost. For heterogeneous DSA, however, computing the action of the off-diagonal blocks requires
a full parallel transport update. Thus, a second option for
heterogeneous DSA preconditioning is a block-diagonal preconditioner, eliminating the need for the
additional update in Algorithm \ref{alg:tri}. Define the preconditioner
\begin{align} \label{eq:approxD}
\mathcal{D}_D^{-1} & = \begin{bmatrix} I & \mathbf{0} \\ \mathbf{0} & I + D_{ff}^{-1}\Sigma_{t,f}\end{bmatrix}.
\end{align}
for which the action can be described as follows.
\begin{algorithm}[Diagonal heterogeneous DSA]\label{alg:diag}
\text{ }
\begin{enumerate}
\item Let $\boldsymbol{\varphi}_f$ denote $\boldsymbol{\varphi}$ restricted to the thick region, and
update $\boldsymbol{\varphi}_f \mapsfrom \boldsymbol{\varphi}_f + D_{ff}^{-1}\Sigma_{t,f}\boldsymbol{\varphi}_f$,
where $D_{ff}$ is the DSA diffusion preconditioner over the thick region.
\end{enumerate}
\end{algorithm}

When preconditioning by inverting the diagonal blocks of a $2\times 2$ operator, block-diagonal
preconditioned minimal-residual iterations generally converge to a given tolerance in roughly twice as
many iterations as block-triangular preconditioners, {due to not capturing off-diagonal coupling
in the preconditioner} \cite{2x2}. For more general approximations (as
used here), the difference in convergence between block-diagonal and block-triangular preconditioning
is more complicated. In particular, there are situations where block-diagonal preconditioning can converge
as fast as block-triangular, or many times slower \cite{diag}. In almost all problems we have tested, 
heterogeneous DSA appears to fall into the former -- block-diagonal preconditioning results in
convergence as good as, or almost as good as, block-triangular preconditioning, at a fraction of the
cost.  Block-triangular preconditioning has resulted in fewer iterations on a few problems where
$\Sigma_{s,s} \sim \mathcal{O}(1)$ is considered ``thin'' (for example, see Figure \ref{fig:cp_comp1}
in Section \ref{sec:results:pipe:2}),
but due to the auxiliary update, the overall time-to-solution remained longer than with block-diagonal
preconditioning. {Block-diagonal preconditioning has the additional advantage over block-triangular
preconditioning that it can be used with conjugate gradient or MINRES if the scalar flux problem
is symmetrized, as in \cite{Gupta:2002ft,Chang:2007go}.}

Thus, the main algorithm proposed here is very simple: only apply DSA preconditioning on the
thick region. It is immediately apparent why this is also compatible with voids, because if
$\Sigma_t = \mathbf{0}$ on a certain region in the problem domain, a diffusion approximation
(typically depending on $\Sigma_t^{-1}$) is neither formed nor inverted on that region. It is
is important to note that heterogeneous DSA is \textit{not the same} as inverting a global
DSA preconditioner and only applying updates to the thin region. 

\begin{remark}[Block-diagonal vs. block-triangular heterogeneous DSA]
Examples were constructed in \cite{diag} where block-diagonal preconditioned minimal-residual
iterations took up to $10\times$ more iterations to converge than with block-triangular preconditioning.
If such a problem were to arise in transport, the auxiliary transport update in block-triangular heterogeneous
DSA would likely be worth the added computational cost. Moreover, block-diagonal preconditioners
typically rely more heavily on Krylov acceleration for convergence. When preconditioning source iteration
without Krylov acceleration, block-triangular preconditioning likely provides a more robust method. 
\end{remark}

\section{Implementation}\label{sec:imp}

\subsection{DSA discretizations}\label{sec:imp:disc}

Conceptually, the heterogeneous preconditioning technique developed in Section \ref{sec:alg} is flexible
in terms of underlying transport and diffusion discretization. This paper focuses on discontinuous Galerkin (DG) 
discretizations of the spatial transport equation and thus, for consistency, DG DSA discretizations as well.
DG is of particular interest because the discontinuous finite-element framework is amenable to traditional
techniques such as sweeping, and allows for high-order accuracy, including discretizing on high-order
curvilinear finite elements.

DSA preconditioning for DG discretizations of S$_N$ transport has been considered in a number of
papers, perhaps originally in \cite{Warsa:2003ug} and more recently considering high-order 
discretizations in \cite{18dsa,Wang:2010dva}. Here, we build our DSA discretization based on the discrete
analysis performed in \cite{18dsa}, and extending modified stabilization ideas from \cite{Wang:2010dva}.
Following the standard derivation of DG discretization, the S$_N$ transport equations in the
diffusive-limit scaling \cite{Larsen:1982hh} can be written in the discrete form
\begin{equation}
\boldsymbol{\Omega}_{d}\cdot\mathbf{G}\boldsymbol{\psi}^{(d)}+F^{(d)}\boldsymbol{\psi}^{(d)}+\frac{1}{\varepsilon} \Sigma_{t}\boldsymbol{\psi}^{(d)}-\frac{1}{4\pi}\left(\frac{1}{\varepsilon} \Sigma_{t}-\varepsilon  \Sigma_{a}\right)\boldsymbol{\varphi}=\frac{1}{4\pi}\left(\varepsilon\boldsymbol{q}_{\text{inc}}^{(d)}+\varepsilon\boldsymbol{q}^{(d)}\right),\label{eq:Tmat1}
\end{equation}
where $\varepsilon$ is the characteristic mean-free path, $\boldsymbol{\psi}_d$ is the angular flux
vector for the $d$th discrete ordinate direction, and $\boldsymbol{\varphi}$ is the scalar flux, given by
$\boldsymbol{\varphi}=\sum_{d}\omega_{d}\boldsymbol{\psi}^{(d)}$, for quadrature weights $\{\omega_d\}$.
{Right-hand side vectors $\boldsymbol{q}_{\text{inc}}^{(d)}$ and $\boldsymbol{q}^{(d)}$ correspond
to the linear forms
\begin{align*}
\left[ \boldsymbol{q}_{\text{inc}}^{(d)} \right]_m & = 
  \sum_{\Gamma\in\mathcal{F}}
  \int_{\Gamma} \boldsymbol{\Omega}_{d} \cdot \mathbf{n}
                v_m \psi_{d,\text{inc}}~dS -
  \frac{1}{2}\sum_{\Gamma\in\mathcal{F}}
  \int_{\Gamma} \left|\boldsymbol{\Omega}_{d}\cdot\mathbf{n} \right|
                v_m \psi_{d,\text{inc}}~dS, \\
\left[ \boldsymbol{q}^{(d)} \right]_m & =
  \sum_{\kappa}\int_{\kappa} v_m q^{(d)} d\mathbf{x},
\end{align*}
for finite-element basis $\{ v_m \}$.}
Here $\Sigma_a, \Sigma_s$, and $\Sigma_t$ are mass matrices corresponding to coefficients $\sigma_a(\mathbf{x})$,
$\sigma_s(\mathbf{x})$, and $\sigma_t(\mathbf{x})$, respectively; $F^{(d)}$ is the DG face matrix based
on upwinding with respect to $\boldsymbol{\Omega}_d$; and
$\boldsymbol{\Omega}_{d}\cdot\mathbf{G}$ is the element-wise discretization of advection in direction
$\boldsymbol{\Omega}_{d}$. For more details, see \cite{18dsa}.

The analysis in \cite{18dsa} proves that a symmetric interior penalty (SIP) discretization provides a
robust DSA preconditioner in the optically thick limit, $\varepsilon \ll 1$, given by
{
\begin{equation}
D_{\textnormal{SIP}} = \frac{1}{\varepsilon}F_{0} + \mathbf{G}^{T}\cdot\mathbf{I}\Sigma_{t}^{-1}\cdot\mathbf{G} +
	\mathbf{F}^T_{1}\cdot  \Sigma_{t}^{-1}\mathbf{G} + \mathbf{G}^T\cdot  \Sigma_{t}^{-1}\mathbf{F}_{1} +  \Sigma_{a},\label{eq:SIP}
\end{equation}
where 
\begin{equation*}
F_{0} = \frac{1}{4\pi}\sum_{d}w_{d}F_{\left\llbracket \right\rrbracket}^{(d)},\,\,\,\,\,
\mathbf{F}_{1}=\frac{1}{4\pi}\sum_{d}w_{d}\boldsymbol{\Omega}_{d}F^{(d)}_{\{\}},\,\,\,\,
\mathbf{I}=\frac{1}{4\pi}\sum_{d}w_{d}\boldsymbol{\Omega}_{d}\boldsymbol{\Omega}_{d},
\end{equation*}
and we perform discrete angular integration of the bilinear forms
\begin{equation*}
\mathbf{v}^T F^{(d)}_{\{\}}\mathbf{u} = - \sum_{\Gamma\in\mathcal{F}}\int_{\Gamma}\mathbf{\Omega}_d\cdot\mathbf{n}\left\llbracket u\right\rrbracket \left\{ v\right\} dS,\,\,\,\,\,
\mathbf{v}^T F_{\left\llbracket \right\rrbracket}^{(d)}\mathbf{u} = \sum_{\Gamma\in\mathcal{F}}\int_{\Gamma}\frac{1}{2}\left|\boldsymbol{\Omega}_{d}\cdot\mathbf{n}\right|\left\llbracket u\right\rrbracket \left\llbracket v\right\rrbracket dS .
\end{equation*}
}

Numerical results in \cite{18dsa} then indicate that dropping the face term
$\mathbf{G}^T\cdot  \Sigma_{t}^{-1}\mathbf{F}_{1}$, leading to a nonsymmetric interior
penalty (NIP) discretization, provides a more robust preconditioner in non-optically thick material,
\begin{equation}
D_{\textnormal{NIP}} = \frac{1}{\varepsilon}F_{0} + \frac{1}{3}\mathbf{G}^{T}\cdot  \Sigma_{t}^{-1}\mathbf{G} 
{
+ \mathbf{F}^T_{1} 
}
\cdot \Sigma_{t}^{-1}\mathbf{G} +  \Sigma_{a} \label{eq:NIP}.
\end{equation}
For linear DG discretizations and straight-edged meshes, $D_{\textnormal{NIP}}$ is analogous to
the Warsa-Wareing-Morel consistent diffusion discretization developed in \cite{Warsa:2002tq}.
Moreover, with some work, one can show that $D_{\textnormal{NIP}}$ can also be derived by
integrating the first two angular moments of \eqref{eq:sntransport}.
We use $D_{\textnormal{NIP}}$ as our baseline DSA discretization here, as we have found it to be
far more robust on heterogeneous domains than $D_{\textnormal{SIP}}$.

\subsubsection{Modified penalty coefficient}\label{sec:imp:imp:mip}

The~interior penalty term $F_0$ properly enforces continuity of the~solution 
in \eqref{eq:SIP} and \eqref{eq:NIP} in the~optically thick limit, 
$\varepsilon \ll 1$. However, the~penalizing term scales as
$\frac{1}{\varepsilon} F_0$ and tends to zero in the optically thin limit, 
$\varepsilon \gg 1$, and the~preconditioners \eqref{eq:SIP} and \eqref{eq:NIP} 
are found to be unstable. In particular, by continuity of eigenvalues as a function 
of matrix entries, as $\varepsilon\to 0$, \eqref{eq:SIP} transitions continuously
from being SPD to being an indefinite operator. This can lead to eigenvalues very
close to zero, which are very difficult to precondition. We introduce the following
modified penalizing term,
\begin{equation}
  \tilde{F}_0 = \frac{1}{8\pi}\sum_{d}w_{d} 
  \sum_{\Gamma} \int_\Gamma 
  \max\left(\frac{c(p^{up})}{\sigma_t^{up} h^{up}}, 1\right)
  |\boldsymbol{\Omega}_{d}\cdot\boldsymbol{n}| 
  \llbracket u \rrbracket \llbracket v \rrbracket dS.
  \label{eq:ModifiedF0}
\end{equation}
Here, the~sum is over mesh faces $\Gamma$, $\boldsymbol{n}$ is 
the~normal vector to the~face $\Gamma$, $h$ is the~cell size, $c(p) = C p (p+1)$
for constant $C$ (we use $C=0.1$ via testing for the fastest convergence), $p$ is
the~polynomial order of test and trial functions $v$ and $u$, and
$\llbracket~\rrbracket$ is the~jump operator. The~label $up$ represents the~value on
the~upwind side at  the~integration point of the~face $\Gamma$ with respect to 
direction $\boldsymbol{\Omega}_d$. Note, that \eqref{eq:ModifiedF0} coincides 
with the~modified interior penalty formulation in \cite{Wang:2010dva}, but 
it introduces a~general formulation valid on curvilinear meshes,
where the~upwind value can vary along a~highly curved face. 
The~modified form of the~penalty term
\eqref{eq:ModifiedF0} cancels the~scaling $\varepsilon^{-1}$ in the~thin limit
of \eqref{eq:SIP} and \eqref{eq:NIP} and the~penalization successfully 
enforces the~continuity for any transport regime. {mNIP is used to denote
the resulting modified nonsymmetric interior penalty discretization.}
{Note, \eqref{eq:ModifiedF0} is not compatible with regions of void, $\sigma_t = 0$.
However, full DSA is incompatible with voids anyways and, as discussed in the following
subsection, a non-modified coefficient is most appropriate for heterogeneous
preconditioning.}

\subsection{Heterogeneous DSA implementation}\label{sec:imp:imp}

Implementing heterogeneous DSA requires identifying and discretizing on a ``thick'' subdomain. In
doing so, we want (i) the method to be amenable to easy addition to existing codes, (ii) the DSA
discretization to be an effective preconditioning on the thick subdomain, and (iii) the DSA
discretization on the subdomain to be solvable using standard AMG techniques.

Building on these ideas, we first define a tolerance, $\eta$, where mesh elements such that
$\sigma_s(\mathbf{x}) \geq \eta$ are considered thick and mesh elements such that
$\sigma_s(\mathbf{x}) < \eta$ are considered thin. On every processor, it is straightforward
to evaluate $\sigma_s$ on each mesh element and track which elements are considered
thick and thin. Next, suppose we discretize a global DSA operator (as if we are doing normal
``full'' DSA) and then extract the principle submatrix corresponding to thick DOFs. The resulting 
matrix represents the thick region throughout the interior of the subdomain, with weakly enforced
Dirichlet boundary conditions. 
{If one writes carefully the stencil of heterogeneous NIP, it can be seen that 
the~principle submatrix uses exclusively $\sigma_s(\mathbf{x}) \geq \eta$ corresponding to 
thick elements. On the other hand, heterogeneous mNIP introduces thin values $\sigma_s(\mathbf{x}) < \eta$
into the thick principle submatrix for all elements laying on the~thick-thin boundary. This modification
at the boundary makes for less robust preconditioning, which suggests (non-modified) NIP as
the appropriate heterogeneous discretization of in practice.}

The \textit{hypre} library \cite{Falgout:2002vu} now supports
parallel extraction of submatrices from a parallel \textit{hypre} matrix, given a (local) list
of row/column-indices. \textit{hypre} is one of the standard parallel multigrid libraries
and, if not already being used for linear solves in a software package, is easily connected to
facilitate this operation.\footnote{Note, for an optimized implementation, it is likely preferable
to only discretize and construct the DSA matrix on the thick subdomain, particularly if the thick
region is very small (for example, see Section \ref{sec:results:hohlraum}). However, for many
finite-element libraries or existing code bases, this ends up being more intrusive to implement.}

\section{Numerical results}\label{sec:results}

This section demonstrates the efficacy of the heterogeneous DSA approach on difficult problems in
transport with heterogeneous cross sections. We emphasize both the superior and more robust preconditioning
that heterogeneous DSA offers over traditional ``full'' DSA, as well as the important property that the
resulting linear systems are more tractable to solve, typically with $\mathcal{O}(1)$ AMG iterations,
independent of problem size.

AMG solves are performed using {BoomerAMG} in the \textit{hypre} library \cite{Falgout:2002vu}. It is
well-known that the choice of parameters is important for AMG. Here we use a V$(1,1)$-cycle with 
$\ell^1$ hybrid Gauss-Seidel relaxation (\textit{hypre} type 8) \cite{Baker:2011ib}, extended$+i$
(distance-two) interpolation (\textit{hypre} type 6) \cite{DeSterck:2008fc}, HMIS coarsening with
one level of aggressive coarsening (\textit{hypre} type 10) \cite{DeSterck:2006et}, and a strength
threshold of $0.05$.
These are similar parameters to the ``optimized'' parameters chosen for PDT transport runs in
\cite{hanus2019weak}, with the important modification of using symmetric $\ell^1$ hybrid Gauss-Seidel
relaxation. In transport simulations, it is often the case that the spatial problem size per processor
is relatively small, and $\ell^1$-relaxation can be faster/more robust \cite{Baker:2011ib}. 
AMG is used as a preconditioner for GMRES, which is solved to $10^{-4}$ relative residual 
tolerance. 

Spatial linear transport equations are discretized using an upwind DG finite element method,
constructed using the MFEM finite element library \cite{mfem-library}, and an
S$_N$ discretization is used in angle. As discussed in Section \ref{sec:alg}, the angular flux
vectors are eliminated from the system and only the scalar flux is stored and iterated on. 
Iterations are accelerated using fGMRES \cite{Saad:1993fc}, an important choice when using
AMG-preconditioned Krylov methods to solve the DSA matrix, as if the residual is not converged
to very small tolerances, each iteration is actually a different preconditioner. It has been noted in multiple
papers, originally in \cite{Warsa:2004gt,Gupta:2002ft}, that Krylov acceleration is important for
heterogeneous domains. The same is true for the new heterogeneous DSA algorithm proposed
here, and we do not present unaccelerated (fixed-point) results. The code does not have a
traditional parallel sweeping implementation, but angular flux problems are solved using the
nonsymmetric AMG method based on approximate ideal restriction (AIR)
\cite{air2,air1} in the \textit{hypre} library \cite{Falgout:2002vu}, as studied in
\cite{hanophy2020}.

\subsection{The crooked pipe problem}\label{sec:results:pipe}

The first problem we consider is the so-called crooked-pipe problem, originally introduced in \cite{gentile2001implicit}
and discussed as a benchmark for DSA in \cite{SmedleyStevenson:2015wa}. This is a steady state test problem
for thermal radiative transport with a single energy group, and purely isotropic scattering throughout the
domain. The domain is surrounded by vacuum, has a uniform isotropic radiation
field, and has an inward isotropic source $10^4$ times stronger than the radiation field. 
In \cite{SmedleyStevenson:2015wa}, the problem is introduced with two scattering cross sections.
Here we modify the domain to have five regions, shown in Figure \ref{fig:pipe}, to allow for a larger variety
of heterogeneities. Scattering cross section $\sigma_s(\mathbf{x})$ is defined to be piecewise constant
over the subdomains shown in Figure \ref{fig:pipe}, and total cross section is then defined as
$\sigma_t(\mathbf{x}) = \sigma_s(\mathbf{x}) + 1/cdt$, where $c$ 
is the~speed of light and $dt$ the~time step of the~transport solution.
We use an S$_8$ angular discretization with {40 angles in the quadrature set (using the symmetry in 2D)}, and
a 2nd-order DG finite element discretization of the linear spatial transport equation for each angle,
unless otherwise noted. 

\begin{figure}[!ht]
    \begin{center}
            \includegraphics[width=0.65\textwidth]{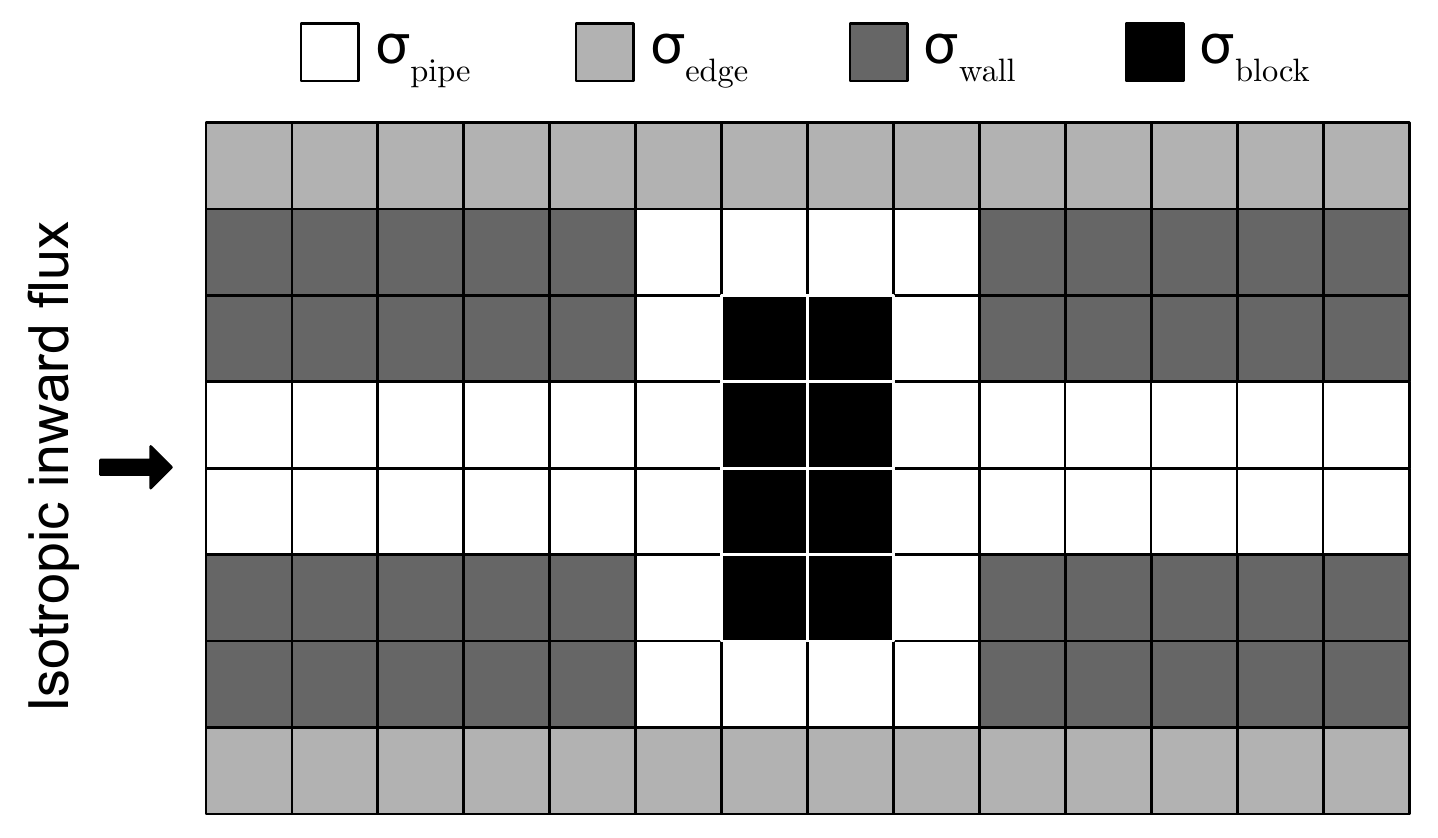}
    \caption{A modification of the crooked pipe problem proposed in \cite{SmedleyStevenson:2015wa}, here allowing
    for $\sigma_s(\mathbf{x})$ to have four different regions, denoted by varying shades from white to black.
    {Each cell in the mesh is size $0.5\times 0.5$.}}
    \label{fig:pipe}
    \end{center}
\end{figure}

\subsubsection{Two cross sections}\label{sec:results:pipe:2}

To start, we consider only two cross sections \cite{gentile2001implicit}, one for the pipe
(consistent with Figure \ref{fig:pipe}) and one for the wall, consisting of everything outside
of the pipe (wall, edge, and block as denoted in Figure \ref{fig:pipe}). Consider the parameters
used in \cite{SmedleyStevenson:2015wa},
\begin{align*}
\sigma_{\cancel{\textnormal{pipe}}} = 200, \hspace{5ex} \sigma_{\textnormal{pipe}} = 0.2, \hspace{5ex}cdt = 1000,
\end{align*}
{where $\sigma_{\cancel{\textnormal{pipe}}}$ corresponds to
all regions outside of the pipe.} Because these are only moderately heterogeneous
coefficients, we instead fix $\sigma_{\cancel{\textnormal{pipe}}} = 200$ and $cdt = 1000$ as above, then scale
$\sigma_{\textnormal{pipe}}$ from zero to 100, to test a wide range of heterogeneities. Numerical tests compare
the newly developed heterogeneous DSA with full DSA, using NIP and mNIP DG diffusion discretizations
(see Section \ref{sec:imp:disc}).

We start by applying {heterogeneous DSA only} to $\sigma_{\cancel{\textnormal{pipe}}}$
{(that is, outside of the pipe)}, even for $\sigma_{\textnormal{pipe}} \geq 1$,
in order to determine what values of $\sigma$ require DSA for good convergence. Figure \ref{fig:dsa}
shows the number of fGMRES iterations to converge transport iterations to a relative residual tolerance of
$10^{-12}$ (with DSA preconditioning applied directly using SuperLU \cite{li1999superlu}), and Figure
\ref{fig:amg} shows the number of AMG iterations to solve a representative
DSA matrix to  $10^{-12}$ relative residual.
}

\begin{figure}[!ht]
    \centering
    \begin{center}
        \begin{subfigure}[t]{0.475\textwidth}
            \includegraphics[width=\textwidth]{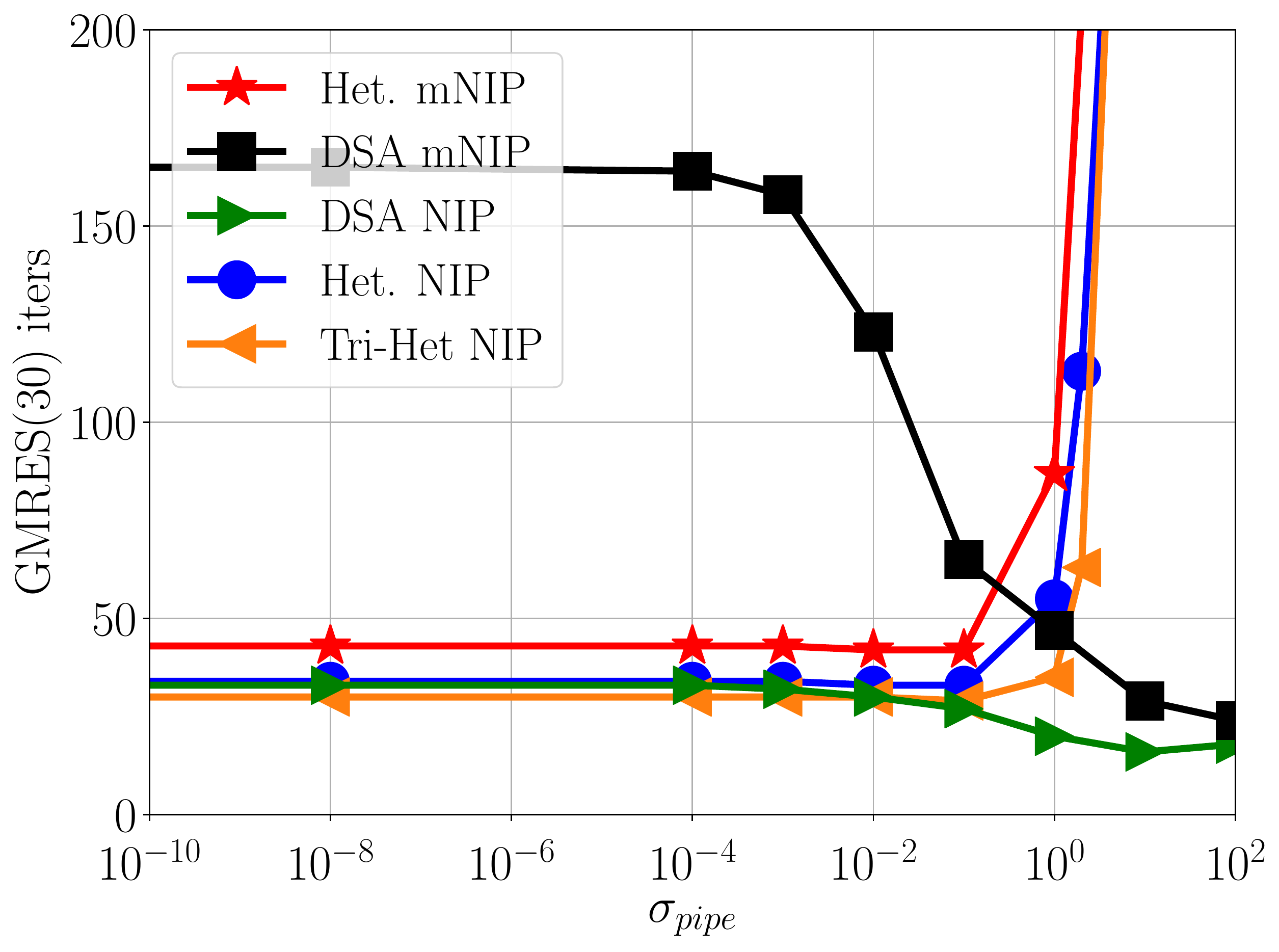}
	\caption{Comparison of fGMRES transport iterations to $10^{-12}$ relative residual
	tolerance for various DSA and heterogeneous DSA preconditioners as a function
	of $\sigma_{\textnormal{pipe}}$. {DSA matrices are inverted directly.}}
	\label{fig:dsa}
        \end{subfigure}\hspace{2ex}
        \begin{subfigure}[t]{0.47\textwidth}
            \includegraphics[width=\textwidth]{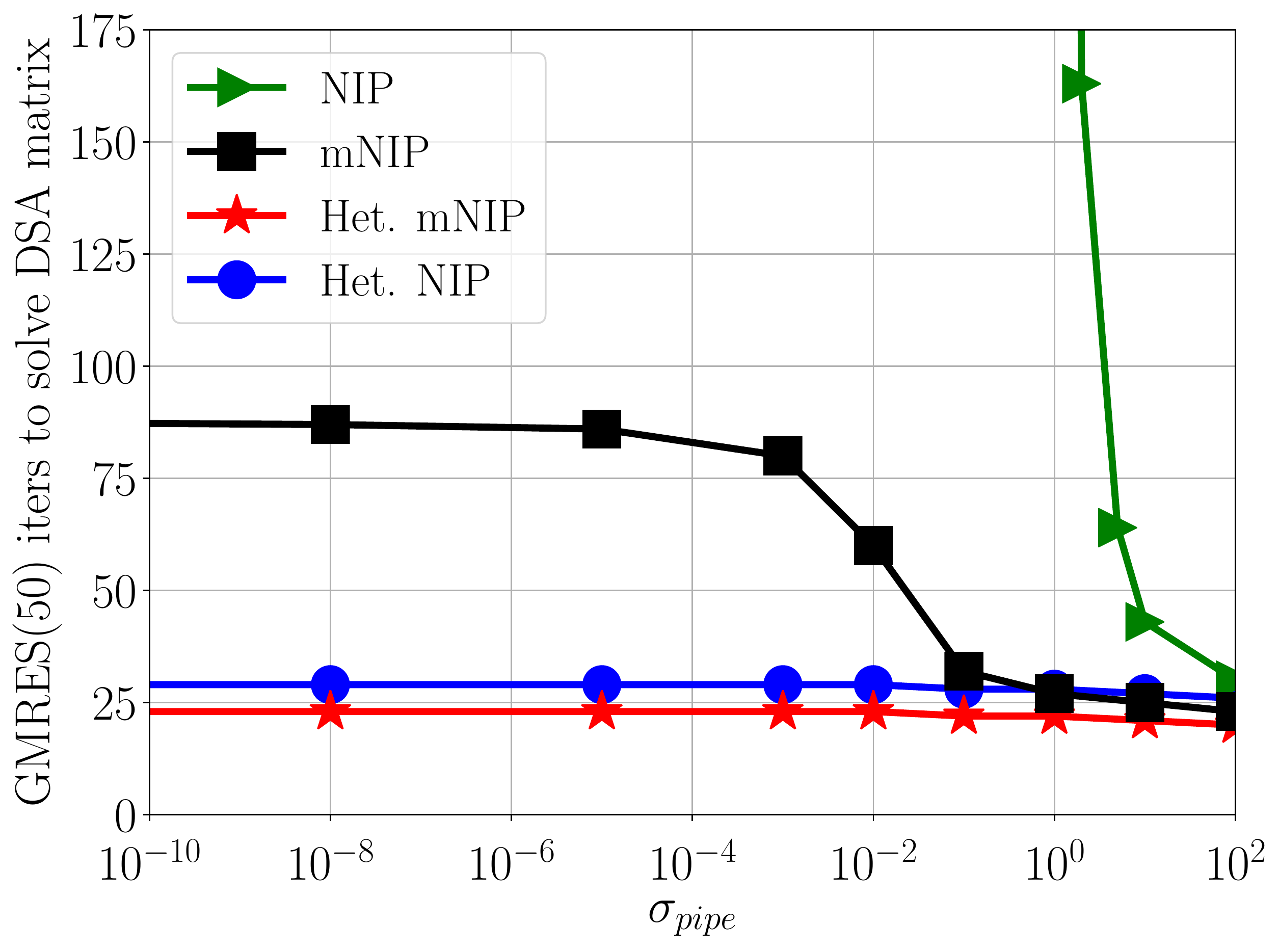}
	\caption{Total number of GMRES iterations, preconditioned by AMG,
	to solve the DSA matrix to a relative residual tolerance of $10^{-12}$, {for a single
  transport iteration}.}
        \label{fig:amg}
        \end{subfigure}
    \end{center}
    \caption{{Results related to DSA} preconditioners applied to the crooked pipe problem with
    $\sigma_{\cancel{\textnormal{pipe}}} = 200$ and $cdt = 1000$, where heterogeneous DSA is only
    applied outside of the pipe. {Note, results in Figures \ref{fig:dsa} and \ref{fig:amg} are
    independent, separately demonstrating convergence of transport iterations and AMG.}}
    \label{fig:cp_comp1}
\end{figure}
\noindent There are a number of {observations} to note from Figure \ref{fig:cp_comp1}:
\begin{itemize}
\item Full DSA with NIP is robust across the entire range of $\sigma_{pipe}$ (green right-pointing
triangles, Figure \ref{fig:dsa}); however, it is very difficult to solve for $\sigma_{pipe}\lessapprox 1$
(green right-pointing triangles, Figure \ref{fig:amg}), and we have tried many different variations of 
AMG and AMG-like techniques. 

\item Full DSA with mNIP is solvable with AMG (black squares, Figure \ref{fig:amg}), {but the
preconditioned transport iteration requires about $5\times$ more iterations for convergence 
on the highly heterogeneous problems (black squares, Figure \ref{fig:dsa}) compared with NIP
preconditioning.} 

\item {Results here are shown for heterogeneous DSA with NIP and mNIP. Recall from Section
\ref{sec:imp:imp} that NIP provides a more natural boundary condition for heterogeneous
preconditioning and will be used for most results. However, it should be noted that regions with small
absorption $\sigma_a \ll 1$ and roughly $\sigma_s\in[1,5]$ require preconditioning of the transport
iteration, but wherein the NIP discretization can be difficult to solve using AMG. Here we demonstrate \
that mNIP may be an effective option as well, although the heterogeneous preconditioning would
be less robust in theory (and, in our experience, also less robust in practice).}

\item Heterogeneous DSA is applied \textit{only outside of the pipe} for these tests, even when
the pipe is moderately thick. {As expected, convergence of heterogeneous DSA deteriorates
rapidly when $\sigma_{pipe}$ is not thin, as such regimes are known to require preconditioning.}
In practice, there is a switch so that DSA can be applied
to all regions such that $\sigma_s(\mathbf{x}) \geq \eta$ for some tolerance $\eta$. If 
$\sigma_{\textnormal{pipe}} > \eta$ for all $\mathbf{x}$, heterogeneous DSA breaks
down to traditional DSA. 
Figure \ref{fig:cp_comp1} suggests a tolerance of $\eta \approx 1.0$ to minimize
the total number of iterations, while ensuring het-NIP DSA is solvable with AMG.
{Note, here we have small absorption, $\sigma_a \sim 1/cdt = 10^{-3}$, and see
comparable results down to zero absorption (e.g., see later results in Table
\ref{tab:crookedpipe1}). For regions with large absorption, a larger $\eta$ may
be possible (resulting in a
smaller thick subdomain), but is also likely unnecessary, that is, $\eta \approx 1$
will provide effective preconditioning.}

\item {For $\sigma_{\textnormal{pipe}} < 0.1$, heterogeneous DSA with NIP yields identical
iteration counts as full DSA with NIP and up to a $5\times$ reduction in iterations
compared to DSA with mNIP. Moreover, AMG solvers are robust when applied to heterogeneous DSA,
{only requiring 20--28 iterations} to obtain 12 digits of accuracy (blue circles and
red stars, Figure \ref{fig:amg}).}

\item Triangular heterogeneous DSA (Tri-Het. DSA) only offers a significant reduction in iteration
count over diagonal heterogeneous DSA for a small range of {$\sigma_{\textnormal{pipe}} \in [1,5]$.
At best, triangular heterogeneous DSA requires $1.8\times$ less iterations than diagonal
heterogeneous DSA for $\sigma_{\textnormal{pipe}}=5$, but with a likely comparable added
computational cost of an additional sweep each iteration. Thus, it seems unlikely
the triangular version will be a better choice in practice, consistent with the discussion in
Section \ref{sec:alg:het}.}

\end{itemize}

Based on these results, let the tolerance $\eta = 1$ moving forward. The following sections
consider harder problems using full DSA with NIP and mNIP and diagonal heterogeneous DSA with NIP.

\subsubsection{Multiple cross sections}\label{sec:results:pipe:3}

Above, we took a benchmark problem and demonstrated the robustness and benefits of heterogeneous DSA.
In this section, we take the same problem but use the full five regions of cross section indicated in Figure \ref{fig:pipe}
to allow for a wider variety of heterogeneities and (hypothetically) more difficult problems. Table \ref{tab:problems}
provides five different sets of scattering cross sections, $\sigma_s(\mathbf{x})$, in the five subregions in
Figure \ref{fig:pipe}. Letting $\eta = 1$ distinguish between thick and thin regions, the percentage of the
domain for each is also given.
{
\begin{table}[!h]
\renewcommand{\tabcolsep}{0.195cm}
\small
\centering
\begin{tabular}{|| c | c c c c | c c ||}\Xhline{1.25pt}
& \multicolumn{4}{c|}{Parameters} & \multicolumn{2}{c||}{Thin/thick} \\\hline
Problem \# & $\sigma_{\textnormal{pipe}}$ & $\sigma_{\textnormal{wall}}$ & $\sigma_{\textnormal{edge}}$ & $\sigma_{\textnormal{block}}$ 
	& \% Thin & \% Thick \\\hline\hline
\#1 & 1e-3 & 500 & 1e-4 & 100 & 55\% & 45\% \\
\#2 & 0.1 & 200 & 200 & 5 & 32\% & 68\% \\
\#3 & 1e-4 & 10 & $500$ & 0.1 & 39\% & 61\% \\ 
\#4 & 1e-6 & 0.1 & $100$ & 1000 & 68\% & 32\% \\
\#5 & 1e-4 & 10 & $ 500 $ & 100 & 32\% & 68\% \\\hline\hline
\Xhline{1.25pt}
\end{tabular} 
\caption{Coefficients $\sigma_s(\mathbf{x})$ for five test problems.}
\label{tab:problems}
\end{table}
}

Results {using 32 cores and $\approx 2500-5000$ DOFs/core} are shown in Table \ref{tab:crookedpipe1}.
Tests are run for $cdt = 1, 1000,$ and $10^8$, as well as the limiting case of zero absorption,
$\sigma_a = 0$, and for 2nd- and 4th-order finite elements. All DSA matrices are
inverted directly using SuperLU \cite{li1999superlu} in order
to compare the preconditioning of S$_N$ transport, independent of the solvability of the DSA matrix. 

{
\begin{table}[!htb]
\renewcommand{\tabcolsep}{0.175cm}
\small
\centering
\begin{tabular}{||c || c c c || c | c || c ||}\Xhline{1.25pt}
& \multicolumn{3}{c||}{Parameters} &  \multicolumn{2}{c||}{Full DSA} & {Het-DSA} \\\hline
& Pr \# & $cdt$ & order & {NIP} & {mNIP} & {NIP}  \\\hline
\textit{1}& \#1 & 1 & 2 & ~19 & ~56 & ~15 \\
\textit{2}& \#2 & 1 &2 & ~18 & ~41 & ~15 \\
\textit{3}& \#3 & 1 & 2 & ~16 & ~41 & ~13 \\
\textit{4}& \#4 & 1 & 2 & ~19 & ~58 & ~15 \\
\textit{5}& \#5 & 1 &2 & ~20 & ~44 & ~15 \\
\hline\hline
\textit{6}& \#1 & $10^3$ & 2& ~31 & 120 & ~25 \\
\textit{7}& \#2 & $10^3$ &2& ~25 & ~61 & ~25  \\
\textit{8}& \#3 & $10^3$ & 2& ~23 & ~96 & ~20 \\
\textit{9}& \#4 & $10^3$ & 2& ~26 & 132 & ~22  \\
\textit{10}& \#5 & $10^3$ &2& 102 & 143 & ~25 \\
\hline\hline
\textit{11}& \#1 & $10^8$ & 2& ~34 & 123 & ~25 \\
\textit{12}& \#2 & $10^8$ &2& ~25 & ~64 & ~25  \\
\textit{13}& \#3 & $10^8$ & 2& ~28 & ~98 & ~21 \\
\textit{14}& \#4 & $10^8$ & 2& 116 & 134 & ~22 \\
\textit{15}& \#5 & $10^8$ & 2& 134 & 147 & ~25 \\
\hline\hline
\textit{16}& \#1 & $\sigma_a = 0$ & 2& ~34 & 153  & ~25 \\
\textit{17}& \#2 & $\sigma_a = 0$ &2& ~25 & ~94  & ~25 \\
\textit{18}& \#3 & $\sigma_a = 0$ & 2& ~24 & 128 & ~21 \\
\textit{19}& \#4 & $\sigma_a = 0$ & 2& 119 & 164 & ~22 \\
\textit{20}& \#5 & $\sigma_a = 0$ & 2& DNC & 177 & ~25 \\
\hline\hline
\textit{21}& \#1 & $\sigma_a = 0$ & 4& ~42 & ~79 & ~25 \\
\textit{22}& \#2 & $\sigma_a = 0$ &4& ~59 & ~40 & ~36  \\
\textit{23}& \#3 & $\sigma_a = 0$ & 4& ~59 & ~57 & ~21  \\
\textit{24}& \#4 & $\sigma_a = 0$ & 4& ~59 & ~73 & ~22 \\
\textit{25}& \#5 & $\sigma_a = 0$ & 4& DNC & ~70 & ~27 \\
\hline\hline
\Xhline{1.25pt}
\end{tabular} 
\caption{{DSA-preconditioned GMRES(30) iterations to $10^{-12}$ relative residual for
Problems \#1 -- \#5, with an S$_8$ angular discretization and 2nd- or 4th-order DG
spatial discretization (DNC denotes did not converge in 500 iterations).}}
\label{tab:crookedpipe1}
\end{table}
}

Note that the heterogeneous DSA method with NIP is as good or better than standard DSA
with NIP or mNIP in all cases. For some problems, heterogeneous DSA takes $5-6\times$ less
iterations (see rows 14/15 and 19) than full DSA, or converges when full DSA with NIP does
not (see rows 20 and 25). {The lack of convergence of DSA with NIP in rows 20 and 25 is likely
due to a (near-)singular NIP DSA matrix that cannot be accurately inverted even with a
direct solver.}

{
\begin{table}[!hbt]
\renewcommand{\tabcolsep}{0.175cm}
\small
\centering
\begin{tabular}{|| c c || c c | c c || c c ||}\Xhline{1.25pt}
\multicolumn{2}{||c||}{Parameters} & \multicolumn{2}{c|}{DSA NIP} & \multicolumn{2}{c||}{DSA mNIP} & \multicolumn{2}{c||}{Het. DSA NIP} \\\hline
Pr \# & $cdt$ & Iters. & AMG it. & Iters. & AMG it. & Iters. & AMG it.\\\hline\hline
\#1 & 1 & 35 & DNC & 39 & 14 & 15 & 7 \\
\#2 & 1 & 27 & DNC & 23 & 8 &  17 & 45 \\
\#3 & 1 & 34 & DNC  & 29 & 10 & 13 & 29 \\
\#4 & 1 & 38 & DNC  & 42 & 10 & 13 & 9 \\
\#5 & 1 & 37 & DNC  & 30 & 10 & 14 & 30 \\
\hline\hline
\#1 & $10^3$ & DNC & DNC  & 73 & 99 & 25 & 9 \\
\#2 & $10^3$ & DNC & DNC  & 40 & 13 & 29 & 180\\
\#3 & $10^3$ & DNC & DNC  & 55 & 33 & 22 & 36 \\
\#4 & $10^3$ & DNC & DNC  & 71 & 26 & 23 & 12 \\
\#5 & $10^3$ & DNC & DNC  & 78 & 34 & 27 & 35 \\
\hline\hline
\#1 & $10^8$ & DNC & DNC  & 71 & 145 & 25 & 9 \\
\#2 & $10^8$ & DNC & DNC  & 40 & 12 & 29 & 180\\
\#3 & $10^8$ & DNC & DNC  & 55 & 42 & 22 & 36 \\
\#4 & $10^8$ & DNC & DNC  & 73 & 43 & 23 & 12 \\
\#5 & $10^8$ & DNC & DNC  & 75 & 41 & 27 & 35 \\
\hline\hline
\Xhline{1.25pt}
\end{tabular} 
\caption{{DSA-preconditioned fGMRES(30) iterations to $10^{-12}$ relative residual
(``Iters.''), and AMG-preconditioned GMRES iterations to apply \textit{a single} DSA preconditioning
to $10^{-4}$ relative-residual (``AMG It.'').}}
\label{tab:crookedpipe2}
\end{table}
}Next we consider a refined spatial mesh, with an S$_8$ angular discretization and 2nd-order 
DG finite elements ($\approx 1$M spatial DOFs), and look at convergence of DSA preconditioned
with NIP, mNIP, and heterogeneous NIP. All {simulations are run on 256 cores}, and DSA
solves are performed using a maximum of 250
AMG iterations. In Table \ref{tab:crookedpipe2}, we see that full DSA with NIP does not converge for
$cdt \gg 1$, largely because 250 AMG iterations provide a very poor approximation to
the diffusion inverse. On the refined mesh, mNIP iterations to convergence actually decrease
(compared with Table \ref{tab:crookedpipe1}) and each iteration only requires a modest number
of AMG iterations. {However, heterogeneous NIP provides faster convergence for the transport
iterations in all cases, typically yielding a $2-3\times$ reduction in iteration count,
for comparable (sometimes less, sometimes more) AMG iterations.}

{It should be noted that AMG struggles on the heterogeneous preconditioning for problem
$\#2$. Interestingly, this is also the easiest problem for full DSA with mNIP. These are
related, and due to the facts that (i) Problem \#2 does not have any particularly thin
regions or particularly big discontinuities in $\sigma$ (see Table \ref{tab:problems}),
where full DSA with mNIP is least effective (see Figure \ref{fig:dsa}); and (ii) Problem \#2
has a region of $\sigma = 5$. Cross-sections $\sigma_s\sim\mathcal{O}(1)$ (with small
$\sigma_a$) are exactly the interface where DSA preconditioning is necessary for good
convergence of source iteration, but AMG applied to (non-modified) DG discretizations
begins to struggle (see Figure \ref{fig:cp_comp1}). It should be noted, applying mNIP
with heterogeneous DSA yields a comparable iteration count to full DSA in terms of
both AMG and source iteration. We do not include these results because for all other tested
problems, heterogeneous DSA with mNIP yields a significant increase in outer transport iteration count 
over heterogeneous DSA with NIP, with only a marginal reduction in AMG iterations.
However, we do note that heterogeneous DSA is best suited for problems with large regions
of $\sigma_s \ll 1$ and/or $\sigma_s \gg 1$.
}

\subsection{Hohlraum}\label{sec:results:hohlraum}
\newcommand{\Cv}{C_v}

The~second problem we consider represents a~realistic setting of the~hohlraum 
chamber with a fuel capsule. This design is used in the~indirect-drive 
approach of {inertial confinement fusion (ICF)}. Radiation transport plays
a~fundamental role in the~indirect-drive scenario, where the~gold wall of the
hohlraum, heated by lasers to a high temperature, radiates x-rays that propagate 
through a~low density Helium fill. This causes a compression of the~fuel capsule 
as the photons are absorbed on its~surface. The capsule is filled with hydrogen, 
contained inside a~thin layer of plastic (CH) on its surface, with a~small hole
on the~right simulating a~filling tube. The~cross section and emissivity of a
material is strongly dependent on its composition and temperature, 
making x-ray radiation transport a highly non-linear process.

\begin{figure}[!ht]
    \begin{center}
        \begin{subfigure}[t]{0.475\textwidth}
	    \centering
            \includegraphics[width=0.95\textwidth]{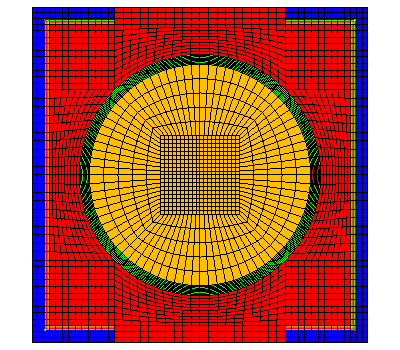}
    \caption{Hohlraum mesh, with materials represented as colors: gold $\leftrightarrow$ blue, 
	helium $\leftrightarrow$ red, CH $\leftrightarrow$ green, and hydrogen $\leftrightarrow$  yellow
  (left-to-right: blue, red, green, yellow, green, red, blue).}
    \label{fig:hohlraum_mesh}
        \end{subfigure}\hspace{2ex}
        \begin{subfigure}[t]{0.47\textwidth}
            \centering
	    \includegraphics[width=0.95\textwidth]{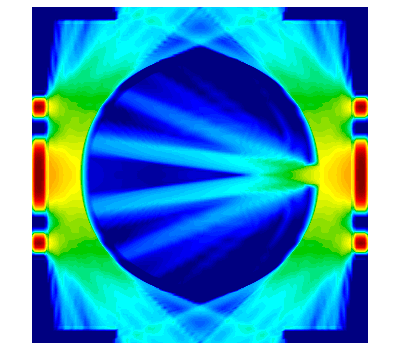}
    \caption{Scalar flux solution (S$_{18}$) due to radiation source from
    the laser-heated gold wall (red regions $\sim$ laser heating).}
    \label{fig:hohlraum_phi}
        \end{subfigure}
    \caption{Indirect-drive ICF approach using a gold cylindrical chamber 
    (hohlraum) filled with helium and a fuel capsule placed in the middle.}
    \label{fig:hohlraum}
    \end{center}
\end{figure}

This can be represented by a~time-dependent thermal radiative transfer problem
modeled by a transport equation with a ``pseudo-scattering'' term
\cite{FleckJCP1971pseudo, LarsenJCP1988gta}, $\sigma_{ps}$,
and a~corresponding source of gray-body radiation, $q_{ps}$:
\begin{align}
  \sigma_{ps} &= \frac{\sigma^2 16\pi a c T^3}
                        {\frac{\Cv}{\Delta t} + 16\pi \sigma a c T^3},
  \label{eq:hohlraum_sigma} \\
  q_{ps} &= \sigma_a a c T^4 - \frac{\sigma^2 16\pi a^2 c^2 T^7} 
          {{\frac{\Cv}{\Delta t} + 16\pi \sigma a c T^3}} ,
  \label{eq:hohlraum_q}
\end{align}
both as a function of coefficient $\sigma$.
Here, $T$ is the plasma temperature, $\Cv$ the~heat capacity, 
$c = 2.9979\times 10^4$, $a = 137.199$ the~radiation constant, and we use 
a~very long time step $\Delta t = 10^{-3} \mu$s, corresponding to $\approx10\%$
of the duration of a~common ICF simulation. The~unit system of 
National Ignition Facility is used \cite{Marinak2008}.
The~hohlraum problem introduces four different materials: gold (blue), 
helium (red), CH (green), and hydrogen (yellow), each with a different 
cross section and heat capacity as shown in Figure \ref{fig:hohlraum_mesh}, 
and with corresponding $\sigma$ and $C_v$ values given in Table \ref{tab:hohlraum_materials}. 
{
\begin{table}[!h]
\renewcommand{\tabcolsep}{0.175cm}
\small
\centering
\begin{tabular}{|| c || c | c | c | c ||}\Xhline{1.25pt}
& hohlraum gold wall & helium fill & capsule CH layer & capsule hydrogen fuel
\\\hline\hline
$\sigma$ & 10$^3$ & 10$^{-3}$ & 10$^2$   & 1 \\
$\Cv$    & 10$^5$ & 10$^{-2}$ & 10$^{3}$ & 1  
\\\hline\hline
\Xhline{1.25pt}
\end{tabular} 
\caption{Opacities, $\sigma_a(\mathbf{x})$, and heat capacities, $\Cv(\mathbf{x})$, for given materials.}
\label{tab:hohlraum_materials}
\end{table}
}

Here we consider a~steady state test of a single time step, modeled by \eqref{eq:sntransport}
with  a~single energy group, and purely isotropic scattering throughout the
domain. In the notation of \eqref{eq:sntransport}, we have
$\sigma_s(\mathbf{x}) = \sigma_{ps}(\mathbf{x})$ and 
$\sigma_t(\mathbf{x}) = \sigma(\mathbf{x})$, with the~radiation source
as $q_d(\mathbf{x}) = q_{ps}(\mathbf{x})$.
The domain is surrounded by vacuum, has an initially 
uniform temperature $T = 10^{-3}$ keV, and the~source $q_{ps}$ on 
the~hohlraum gold wall corresponds to $T = 0.4$ keV. 
The~solution of the scalar flux radiation field produced by a~laser-heated 
gold wall is shown in Figure \ref{fig:hohlraum_phi}.
In large-scale problems the ``pseudo scattering'' equation is solved iteratively, 
and standard iterative methods converge very slowly when either the~time
step or $\sigma$ is large. Extreme heterogeneities are encountered in such cases,
in particular on the~gold-helium and CH-helium material interfaces. 

Table \ref{tab:hohlraum} shows results for DSA and heterogeneous DSA preconditioning of a 4th-order DG
discretization in space, and S$_4$, S$_8$, and S$_{12}$ angular discretizations. {Simulations
are run on 16, 64, and 256 cores, corresponding to the three levels of refinement shown, with
$\approx 5000$ spatial DOFs/core.} 
Performance of AMG on the hohlraum problem is particularly poor for full DSA. Results for full DSA
with NIP are not shown in Table \ref{tab:crookedpipe2} because AMG iterations make no significant
reduction in residual for the DSA solve and, thus, the larger transport iterations do not converge.
AMG is also unable to solve the full DSA with mNIP matrices (1,000 AMG-preconditioned
GMRES iterations reduced the residual less than two orders of magnitude), but apparently 
enough preconditioning is achieved in 250 iterations for the larger transport iterations to
achieve reasonable convergence. However, we would not consider this a robust method, as
results for full DSA with NIP here and in Table \ref{tab:crookedpipe2} indicate that if the AMG iterations
are ineffective, the preconditioning may also be ineffective, and the transport iterations may
not converge. Moreover, it is likely AMG convergence will continue to degrade as the spatial 
problem is further refined, which may eventually prevent convergence of the larger transport
iterations. Note, we also tried multiple scaling constants for mNIP and saw no notable 
improvement in AMG convergence. 

{
\begin{table}[!hbt]
\renewcommand{\tabcolsep}{0.175cm}
\small
\centering
\begin{tabular}{|| c c || c c | c c c||}\Xhline{1.25pt}
& & \multicolumn{2}{c|}{DSA mNIP} & \multicolumn{3}{c||}{Het. DSA NIP} \\\hline
S$_N$ & DOFs & Iters. & AMG it. & Iters. & AMG it. & \% Thick DOFs\\\hline\hline
4 & ~~~78,900 & 52 & DNC & 20 & 6 & 2.3\% \\
4 &  ~~315,600 & 40 & DNC & 23 & 8 & 2.9\% \\
4 & 1,262,400 & 28 & DNC & 20 & 10 & 2.5\%\\ \hline
8 & ~~~78,900 & 52 & DNC & 20 & 6 & 2.3\% \\
8 &  ~~315,600 & 40 & DNC & 23 & 7 & 2.9\%  \\
8 & 1,262,400 & 29 & DNC & 21 & 10 & 2.5\% \\\hline
12 & ~~~78,900 & 52 & DNC & 20 & 6 & 2.3\% \\
12 &  ~~315,600 & 39 & DNC & 23 & 7 & 2.9\%\\
12 & 1,262,400 & 29 & DNC & 21 & 10  &2.5\%\\
\hline\hline
\Xhline{1.25pt}
\end{tabular} 
\caption{\color{black}DSA-preconditioned fGMRES iterations to $10^{-12}$ relative residual, 
AMG-preconditioned fGMRES iterations to $10^{-4}$ relative residual tolerance (with a maximum
250 AMG iterations), and percentage of DOFs marked ``thick,'' for the hohlraum problem.
}
\label{tab:hohlraum}
\end{table}
}

{In contrast, heterogeneous DSA proves to be fast and robust. Heterogeneous DSA converges in
fewer iterations than full DSA for all tested problem sizes and S$_N$ orders, reducing total
iteration counts between 30\% to more than $2.5\times$ for the coarsest mesh. Moreover, AMG
proves effective as a solver for heterogeneous NIP, requiring at most 10 AMG iterations to reduce the
residual four orders of magnitude.} This highlights how interesting the dynamics of the problem
and preconditioning can be: without DSA preconditioning, source iteration on this problem converges
extremely slow. Preconditioning a very small subdomain (less than 3\% of the mesh elements!)
with 10 AMG iterations results in rapid convergence, independent of spatial mesh refinement
or S$_N$ order. 

\section{Conclusions}\label{sec:conc}

This paper introduces a new DSA-like technique to precondition transport iteration
in  highly heterogeneous domains, which is trivially compatible with voids. 
The preconditioning is based on a linear algebraic analysis
rather than the underlying physics, but proves to be at least as fast as standard ``full'' DSA on
all problems we tested, and reduces the iteration count by $5-6\times$ on some examples. 
Moreover, even for robust DSA discretizations based on integrating angular moments,
the resulting linear systems are {typically solvable using $\mathcal{O}(1)$ AMG iterations,
while the same discretization applied to the entire domain (full DSA) can lead to matrices which
are very difficult or intractable to solve.}

Even with optimized AMG parameters and
a state-of-the-art parallel AMG library, the application of DSA can take a remarkably large
portion of the solve phase for S$_N$ transport, as high as 60\% with a small number of
mesh elements per core \cite{hanus2019weak}. For difficult problems like the hohlraum
discussed in Section \ref{sec:results:hohlraum}, AMG simply does not converge, and it
is unclear if the preconditioned transport iterations will converge. In addition to making the
DSA solves tractable, only having to solve a DSA discretization on the moderately thick to
thick subdomains can also significantly reduce the percentage of time the DSA solves
take. The hohlraum problem studied in Section \ref{sec:results:hohlraum} only requires
DSA preconditioning on less than $ 3\%$ of the mesh elements for rapid convergence. 
Finally, the proposed approach is relatively non-intrusive and easy to add to existing
high-performance transport codes. 

Numerical results in Section \ref{sec:results} indicate rather different behavior on different
heterogeneous configurations of $\sigma_t$, in terms of DSA preconditioning and AMG
convergence. Analyzing this problem is nontrivial as it is based on an approximation
to a sum of Schur complement inverses \eqref{eq:approxS0}.
A unified framework for DSA preconditioning, including the heterogeneous approach
developed here as well as a better understanding of how discretization properties and
parameters affect both DSA preconditioning and the performance of multigrid-like solvers,
remains an outstanding issue and long-term objective.

\section*{Acknowledgments}
This material is based upon work supported by the Department of Energy, National Nuclear
Security Administration, under Award Number DE-NA0002376. This work was performed under
the auspices of the U.S. Department of Energy by Lawrence Livermore National Laboratory under
contract DE-AC52-07NA27344 (LLNL-JRNL-802400).
Disclaimer: This document was prepared as an account of work sponsored by an agency of the United States government. Neither the United States government nor Lawrence Livermore National Security, LLC, nor any of their employees makes any warranty, expressed or implied, or assumes any legal liability or responsibility for the accuracy, completeness, or usefulness of any information, apparatus, product, or process disclosed, or represents that its use would not infringe privately owned rights. Reference herein to any specific commercial product, process, or service by trade name, trademark, manufacturer, or otherwise does not necessarily constitute or imply its endorsement, recommendation, or favoring by the United States government or Lawrence Livermore National Security, LLC. The views and opinions of authors expressed herein do not necessarily state or reflect those of the United States government or Lawrence Livermore National Security, LLC, and shall not be used for advertising or product endorsement purposes.

\bibliographystyle{unsrt}
\bibliography{refs.bib}

\end{document}